\documentclass[12pt,twoside]{article}
\usepackage{amsmath, amssymb, amsthm}
\font\fourteenb=cmb10 at 14pt 
\begin{document}
\begin{center}

\vspace*{-1.0cm}\noindent \
 \textbf{MATHEMATICS AND EDUCATION IN MATHEMATICS, 1981}\\[-0.0mm]\
\ Proceedings of the 10th Fourth Spring Conference of

the Union of Bulgarian Mathematicians,\\[-0.0mm]
\textit{ Sunny Beach , April   6 - 9, 1981}\\[-0.0mm]
\font\fourteenb=cmb10 at 12pt \ \

   {\bf \LARGE Inequalities for analytic functions with the derivative in  $H^{1} $
   \\ \ \\ \large Peyo Stoilov}
\end{center}

\


\begin{abstract}
It is proved that if  ${\kern 1pt} {\kern 1pt} f'\in H^{1} $ ,
then
$$\int \limits _{T}\frac{\left|f(\zeta \eta )-f(\overline{\zeta }\eta )\right|}{\left|1-\zeta \right|}  dm(\zeta ){\kern 1pt}={\kern 1pt} \frac{1}{\pi } \int \limits _{0}^{\pi }\frac{\left|f(e^{i(\theta +t)} )-f(e^{i(\theta -t)} )\right|}{2\sin (t\backslash 2)}{\kern 1pt} \le {\kern 1pt} \pi \left\| f'\right\| _{H^{1} } ,{\kern 1pt} {\kern 1pt} {\kern 1pt}{\kern 1pt} \eta =e^{i\theta }.$$
\end{abstract}

     \section{Introduction}

     Let  $A$  be the class of all functions analytic in the unit disc  $D=\left\{\zeta :{\kern 1pt} {\kern 1pt} {\kern 1pt} {\kern 1pt} {\kern 1pt} {\kern 1pt} \left|\zeta \right|{\kern 1pt} <1\right\},$  $m(\zeta )$ - normalized  Lebesgue measure on the circle  $T=\left\{\zeta :{\kern 1pt} {\kern 1pt} {\kern 1pt} {\kern 1pt} {\kern 1pt} {\kern 1pt} \left|\zeta \right|{\kern 1pt} =1\right\}$ . Let  $H^{p{\kern 1pt} {\kern 1pt} {\kern 1pt} {\kern 1pt} } (1\le p\le \infty )$  is the
space of  all functions analytic in  $D$  and satisfying

     $$\left\| f\right\| _{H^{p} } =\mathop{\sup }\limits_{0\le {\kern 1pt} {\kern 1pt} {\kern 1pt} r<{\kern 1pt} {\kern 1pt} {\kern 1pt} 1} \left(\int _{T}\left|f(r\zeta )\right| ^{p} dm(\zeta )\right)^{1/p} <\infty ,{\kern 1pt} {\kern 1pt} {\kern 1pt} {\kern 1pt} {\kern 1pt} {\kern 1pt} {\kern 1pt} {\kern 1pt} {\kern 1pt} {\kern 1pt} {\kern 1pt} {\kern 1pt} {\kern 1pt} {\kern 1pt} {\kern 1pt} {\kern 1pt} {\kern 1pt} {\kern 1pt} {\kern 1pt} {\kern 1pt} {\kern 1pt} {\kern 1pt} {\kern 1pt} {\kern 1pt} 1\le p<\infty ,$$

     $$\left\| f\right\| _{H^{\infty } } =\mathop{\sup }\limits_{z\in D} {\kern 1pt} {\kern 1pt} {\kern 1pt} {\kern 1pt} \left|f(z)\right|<\infty ,{\kern 1pt} {\kern 1pt} {\kern 1pt} {\kern 1pt} {\kern 1pt} {\kern 1pt} {\kern 1pt} {\kern 1pt} {\kern 1pt} {\kern 1pt} {\kern 1pt} {\kern 1pt} {\kern 1pt} {\kern 1pt} {\kern 1pt} {\kern 1pt} {\kern 1pt} {\kern 1pt} {\kern 1pt} {\kern 1pt} {\kern 1pt} {\kern 1pt} {\kern 1pt} {\kern 1pt} p=\infty .$$

     For    $f\in A,{\kern 1pt} {\kern 1pt} {\kern 1pt} {\kern 1pt} {\kern 1pt} {\kern 1pt} f'\in H^{1} $  it follows the Hardy inequality [1, 104-105]:

     $$\sum \limits _{k\ge 1}\left|\hat{f}(k)\right| {\kern 1pt} {\kern 1pt} {\kern 1pt} {\kern 1pt} {\kern 1pt} \le {\kern 1pt} {\kern 1pt} {\kern 1pt} {\kern 1pt} \pi \left\| f'\right\| _{H^{1} }. $$

     In this paper we shall prove an inequality -  integrated  analogue of the  Hardy  inequality and as the application  we shall give simplified proof of the theorem of  S. A. Vinogradov for the bounded Toeplitz  operators on   $H^{\infty } $ [2].
\section{Main results}

\

     {\bf Theorem 1. }{\it  If }$f\in A{\kern 1pt} $ {\it  and  } ${\kern 1pt} f'\in H^{1} $ {\it   then}

 $$\int _{T} \frac{\left|f(\zeta \eta )-f(\overline{\zeta }\eta )\right|}{\left|1-\zeta \right|} dm(\zeta ){\kern 1pt} {\kern 1pt} {\kern 1pt} {\kern 1pt} \le {\kern 1pt} {\kern 1pt} \pi \left\| f'\right\| _{H^{1} } {\kern 1pt} {\kern 1pt} ,{\kern 1pt} {\kern 1pt} {\kern 1pt} {\kern 1pt} {\kern 1pt} {\kern 1pt} {\kern 1pt} {\kern 1pt}{\kern 1pt} \eta \in T.$$

     {\it Proof. }Let  ${\kern 1pt} {\kern 1pt} {\kern 1pt} {\kern 1pt} {\kern 1pt} {\kern 1pt} f'\in H^{1},{\kern 1pt} {\kern 1pt}  {\kern 1pt} {\kern 1pt} {\kern 1pt} {\kern 1pt} \eta \in T.{\kern 1pt} {\kern 1pt} $ Then
\begin{equation}
\int _{T} \frac{\left|f(\zeta \eta )-f(\overline{\zeta }\eta
)\right|}{\left|1-\zeta \right|} dm(\zeta ){\kern 1pt} {\kern 1pt}
{\kern 1pt} {\kern 1pt} \le {\kern 1pt} {\kern 1pt} {\kern 1pt}
\mathop{\lim }\limits_{r\overline{\to 1-}0} {\kern 1pt} {\kern
1pt} \int _{T} \frac{\left|f(r\zeta \eta )-f(r\overline{\zeta
}\eta )\right|}{\left|1-\zeta \right|} dm(\zeta ){\kern 1pt}
{\kern 1pt}.
\end{equation}
     Since ${\kern 1pt} {\kern 1pt} {\kern 1pt} {\kern 1pt} {\kern 1pt} {\kern 1pt} f'\in H^{1} ,$  integrating in parts and applying the Cauchy theorem, we obtain

     $$\begin{array}{l} {\frac{1}{2\pi i} \int _{T} {\kern 1pt} {\kern 1pt} f'(\xi ){\kern 1pt} {\kern 1pt} {\kern 1pt} ln\left|1-\xi \overline{z}\right|^{2} d\xi {\kern 1pt} {\kern 1pt} {\kern 1pt} {\kern 1pt} {\kern 1pt} {\kern 1pt} {\kern 1pt} {\kern 1pt} ={\kern 1pt} } \\ {{\kern 1pt} {\kern 1pt} {\kern 1pt} {\kern 1pt} {\kern 1pt} {\kern 1pt} {\kern 1pt} } \\ {\frac{1}{2\pi i} \int _{T} {\kern 1pt} {\kern 1pt} f'(\xi ){\kern 1pt} {\kern 1pt} {\kern 1pt} ln({\kern 1pt} {\kern 1pt} {\kern 1pt} 1-\xi \overline{z}{\kern 1pt} {\kern 1pt} {\kern 1pt} {\kern 1pt} )d\xi {\kern 1pt} {\kern 1pt} {\kern 1pt} {\kern 1pt} {\kern 1pt} +{\kern 1pt} \frac{1}{2\pi i} \int _{T} {\kern 1pt} {\kern 1pt} f'(\xi ){\kern 1pt} {\kern 1pt} {\kern 1pt} ln({\kern 1pt} {\kern 1pt} {\kern 1pt} 1-\overline{\xi }z{\kern 1pt} {\kern 1pt} {\kern 1pt} {\kern 1pt} )d\xi =} \\ {} \\ {\frac{1}{2\pi i} \int _{T} {\kern 1pt} {\kern 1pt} f'(\xi ){\kern 1pt} {\kern 1pt} {\kern 1pt} ln({\kern 1pt} {\kern 1pt} {\kern 1pt} 1-\xi \overline{z}{\kern 1pt} {\kern 1pt} {\kern 1pt} {\kern 1pt} )d\xi +\frac{1}{2\pi i} \int _{T} {\kern 1pt} {\kern 1pt} f'(\xi ){\kern 1pt} {\kern 1pt} {\kern 1pt} ln({\kern 1pt} {\kern 1pt} {\kern 1pt} \xi -z{\kern 1pt} {\kern 1pt} {\kern 1pt} {\kern 1pt} )d\xi -\frac{1}{2\pi i} \int _{T} {\kern 1pt} {\kern 1pt} f'(\xi ){\kern 1pt} {\kern 1pt} {\kern 1pt} ln{\kern 1pt} {\kern 1pt} \xi {\kern 1pt} {\kern 1pt} d\xi =} \\ {} \\ {\frac{1}{2\pi i} \int _{T} {\kern 1pt} {\kern 1pt} f(\xi ){\kern 1pt} ({\kern 1pt} {\kern 1pt} {\kern 1pt} {\kern 1pt} {\kern 1pt} \frac{\overline{z}}{{\kern 1pt} 1-\xi \overline{z}} {\kern 1pt} {\kern 1pt} {\kern 1pt} {\kern 1pt} {\kern 1pt} -{\kern 1pt} {\kern 1pt} {\kern 1pt} \frac{1}{{\kern 1pt} \xi -z} {\kern 1pt} {\kern 1pt} {\kern 1pt} {\kern 1pt} {\kern 1pt} {\kern 1pt} +{\kern 1pt} {\kern 1pt} {\kern 1pt} {\kern 1pt} \frac{1}{\xi } {\kern 1pt} {\kern 1pt} {\kern 1pt} {\kern 1pt} {\kern 1pt} {\kern 1pt} ){\kern 1pt} {\kern 1pt} {\kern 1pt} d\xi {\kern 1pt} {\kern 1pt} {\kern 1pt} {\kern 1pt} {\kern 1pt} {\kern 1pt} ={\kern 1pt} {\kern 1pt} {\kern 1pt} {\kern 1pt} {\kern 1pt} f(0)-f(z){\kern 1pt} {\kern 1pt} {\kern 1pt} {\kern 1pt} {\kern 1pt} {\kern 1pt} {\kern 1pt} {\kern 1pt} {\kern 1pt} {\kern 1pt} {\kern 1pt} {\kern 1pt} {\kern 1pt} \Rightarrow } \end{array}$$

$$f(z)=f(0)-\frac{1}{2\pi i} \int _{T} {\kern 1pt} {\kern 1pt} f'(\xi ){\kern 1pt} {\kern 1pt} {\kern 1pt} ln\left|1-\xi \overline{z}\right|^{2} d\xi .$$

      $\begin{array}{l} {{\kern 1pt} } \\ {} \end{array}$

      Using last equality  for   $f(z\eta )$ :

     $$f(z\eta )=f(0)-\frac{1}{2\pi i} \int _{T} {\kern 1pt} {\kern 1pt} f'(\xi \eta ){\kern 1pt} {\kern 1pt} {\kern 1pt} ln\left|1-\xi \overline{z}\right|^{2} d\xi ,$$

     we  shall have

     $$\left|f(r\zeta \eta )-f(r\overline{\zeta }\eta )\right|=\left|\frac{1}{2\pi i} \int _{T} {\kern 1pt} {\kern 1pt} f'(\xi \eta ){\kern 1pt} {\kern 1pt} {\kern 1pt} ln\left|\frac{1-\xi r\zeta }{1-\xi r\overline{\zeta }} \right|^{2} d\xi \right|\le $$

     $$\le \int _{T} {\kern 1pt} {\kern 1pt} \left|f'(\xi \eta )\right|{\kern 1pt} {\kern 1pt} {\kern 1pt} \left|ln\left|\frac{1-\xi r\zeta }{1-\xi r\overline{\zeta }} \right|^{2} \right|dm(\xi )=$$

     $$=\int _{T} {\kern 1pt} {\kern 1pt} \left|f'(\xi \eta )\right|{\kern 1pt} {\kern 1pt} {\kern 1pt} \left|ln\frac{(1-r)^{2} +r\left|1-\xi \zeta \right|^{2} }{(1-r)^{2} +r\left|1-\xi \overline{\zeta }\right|^{2} } \right|dm(\xi )\le $$

     $$\le \int _{T} {\kern 1pt} {\kern 1pt} \left|f'(\xi \eta )\right|{\kern 1pt} {\kern 1pt} {\kern 1pt} \left|ln\left|\frac{1-\xi \zeta }{1-\xi \overline{\zeta }} \right|^{2} \right|dm(\xi ).$$

     Further from (1) follows

$$\int _{T} \frac{\left|f(\zeta \eta )-f(\overline{\zeta }\eta )\right|}{\left|1-\zeta \right|} dm(\zeta ){\kern 1pt} {\kern 1pt} \le 2{\kern 1pt} \int _{T}\int _{T}\frac{\left|f'(\xi \eta )\right|}{\left|1-\zeta \right|}{\kern 1pt} \left|ln\left|\frac{1-\xi \zeta }{1-\xi \overline{\zeta }} \right|\right|  {\kern 1pt} dm(\xi )dm(\zeta ){\kern 1pt} \le $$

     $${\kern 1pt} {\kern 1pt} {\kern 1pt} {\kern 1pt} {\kern 1pt} {\kern 1pt} {\kern 1pt} \le \left\| f'\right\| _{H^{1} } \sup \left\{2\int _{T}\left|ln\left|\frac{1-\xi \zeta }{1-\xi \overline{\zeta }} \right|\right| \frac{dm(\zeta )}{\left|1-\zeta \right|} {\kern 1pt} {\kern 1pt} {\kern 1pt} {\kern 1pt} {\kern 1pt} {\kern 1pt} {\kern 1pt} :{\kern 1pt} {\kern 1pt} {\kern 1pt} {\kern 1pt} {\kern 1pt} {\kern 1pt} {\kern 1pt} {\kern 1pt} {\kern 1pt} {\kern 1pt} {\kern 1pt} {\kern 1pt} \xi \in T{\kern 1pt} {\kern 1pt} {\kern 1pt} {\kern 1pt} {\kern 1pt} \right\}{\kern 1pt} .$$

     For end the proof is necessary only to estimate  the integral{\it  }

     $${\kern 1pt} {\kern 1pt} {\kern 1pt} {\kern 1pt} I(\xi )=2{\kern 1pt} {\kern 1pt} \int _{T}\left|ln\left|\frac{1-\xi \zeta }{1-\xi \overline{\zeta }} \right|\right| \frac{dm(\zeta )}{\left|1-\zeta \right|} {\kern 1pt} {\kern 1pt} ,{\kern 1pt} {\kern 1pt} {\kern 1pt} {\kern 1pt} {\kern 1pt} {\kern 1pt} {\kern 1pt} {\kern 1pt} {\kern 1pt} {\kern 1pt} {\kern 1pt} {\kern 1pt} \xi \in T{\kern 1pt} .$$

     We have

$$ I(e^{i\theta } )=2{\kern 1pt} .\frac{1}{2\pi } {\kern 1pt} \int \limits _{-\pi }^{\pi }\left|ln\left|\frac{1-e^{i(\theta +t)} }{1-e^{i(\theta -t)} } \right|\right| \frac{dt}{\left|1-e^{it} \right|} {\kern 1pt} {\kern 1pt} {\kern 1pt} {\kern 1pt} {\kern 1pt} {\kern 1pt} {\kern 1pt} {\kern 1pt} ={\kern 1pt} $$

$$ =\frac{1}{\pi } {\kern 1pt} \int \limits _{0}^{\pi }\left|ln\left|\frac{\sin ((\theta +t)/2)}{\sin ((\theta -t)/2)} \right|\right| \frac{dt}{\sin (t/2)} {\kern 1pt} {\kern 1pt} {\kern 1pt} {\kern 1pt} {\kern 1pt} {\kern 1pt} {\kern 1pt} {\kern 1pt} ={\kern 1pt} $$

$$ =\frac{1}{\pi } {\kern 1pt} \int \limits _{0}^{\pi }\left|ln\left|\frac{tg(\theta /2)+tg(t/2)}{tg(\theta /2)-tg(t/2)} \right|\right| \frac{dt}{\sin (t/2)} {\kern 1pt}=\frac{2}{\pi }{\kern 1pt}\int \limits _{0}^{\infty }\left|ln\left|\frac{tg(\theta /2)+y}{tg(\theta /2)-y} \right|\right| \frac{dy}{y\sqrt{1+y^{2} } }{\kern 1pt} \le $$

$$ \frac{2}{\pi } {\kern 1pt} \int \limits _{0}^{\infty }\left|ln\left|\frac{tg(\theta /2)+y}{tg(\theta /2)-y} \right|\right| \frac{dy}{y} {\kern 1pt} \le {\kern 1pt} {\kern 1pt} \frac{2}{\pi } {\kern 1pt} \int \limits _{0}^{\infty }\left|ln\left|\frac{1+x}{1-x} \right|\right| \frac{dx}{x}{\kern 1pt} ={\kern 1pt} \frac{2}{\pi } \int \limits _{0}^{1}\ln \left(\frac{1+x}{1-x} \right) {\kern 1pt} \frac{dx}{x} {\kern 1pt} ={\kern 1pt} \pi .$$

     {\bf Application of }{\bf Theorem 1.}

     For  $f\in H^{1} $  we denote by  $T_{f} $  the Toeplitz operator on   $H^{\infty } $ ,  defined by

     $$T_{f} h=\int _{T} \frac{\overline{f}(\zeta )h(\zeta )}{1-\overline{\zeta }z} dm(\zeta ),{\kern 1pt} {\kern 1pt} {\kern 1pt} {\kern 1pt} {\kern 1pt} {\kern 1pt} {\kern 1pt} {\kern 1pt} {\kern 1pt} {\kern 1pt} {\kern 1pt} {\kern 1pt} {\kern 1pt} {\kern 1pt} {\kern 1pt} {\kern 1pt} {\kern 1pt} {\kern 1pt} {\kern 1pt} h\in H^{\infty } .$$

     In [2] S. A. Vinogradov   proves   that  if   ${\kern 1pt} {\kern 1pt} {\kern 1pt} {\kern 1pt} f'\in H^{1} $ {\it ,  }then the  Toeplitz  operator{\it   } $T_{f} $   is bounded on  $H^{\infty } $ .

     As the application of   Theorem 1 we shall give a simplified  proof  of  the theorem  of   S. A. Vinogradov  and  estimation  for   $\left\| T_{f} \right\| _{H^{\infty } } .$

     {\bf Theorem 2. }{\it  If } $f\in A{\kern 1pt} $ {\it  and  } ${\kern 1pt} f'\in H^{1} $ {\it   then}

     $$\left\| T_{f} \right\| _{H^{\infty } } \le \left\| f\right\| _{H^{\infty } } +{\kern 1pt} \pi \left\| f'\right\| _{H^{1} } .$$

      {\it Proof.}

     $$\left\| T_{f} \right\| _{H^{\infty } } =\sup \left\{\mathop{\lim }\limits_{r\to 1-0} \left|\int _{T} \frac{\overline{f}(\zeta )h(\zeta )}{1-\overline{\zeta }r\eta } dm(\zeta )\right|:\eta \in T,\left\| h\right\| _{H^{\infty } } \le 1\right\}=$$

     $$=\sup \left\{\mathop{\lim }\limits_{r\to 1-0} \left|\int _{T} \frac{\overline{f}(\zeta \eta )h(\zeta \eta )}{1-r\overline{\zeta }} dm(\zeta )\right|:\eta \in T,\left\| h\right\| _{H^{\infty } } \le 1\right\}\le $$

     $$\le \sup \left\{\mathop{\lim }\limits_{r\to 1-0} \left|\int _{T} \frac{\overline{f}(\zeta \eta )-\overline{f}(\overline{\zeta }\eta )}{1-r\overline{\zeta }} h(\zeta \eta )dm(\zeta )\right|:\eta \in T,\left\| h\right\| _{H^{\infty } } \le 1\right\}+\left\| f\right\| _{H^{\infty } } .$$

     We used, that  $g(z)=\overline{f}(\overline{z}\eta )\in H^{\infty } $  and

     $$\left|\int _{T} \frac{\overline{f}(\overline{\zeta }\eta )h(\zeta \eta )}{1-r\overline{\zeta }} dm(\zeta )\right|\le \left\| f\right\| _{H^{\infty } } \left\| h\right\| _{H^{\infty } } .$$

\

     Further from Theorem 1 follows

   $$\left\| T_{f} \right\| _{H^{\infty } } \le \sup \left\{\left|\int _{T} \frac{\left|f(\zeta \eta )-f(\overline{\zeta }\eta )\right|}{\left|1-\zeta \right|} dm(\varsigma )\right|:\eta \in T\right\}+\left\| f\right\| _{H^{\infty } } \le $$

   $$\le {\kern 1pt} {\kern 1pt} {\kern 1pt} {\kern 1pt} {\kern 1pt} {\kern 1pt} \pi \left\| f'\right\| _{H^{1} } +\left\| f\right\| _{H^{\infty } } .$$

\end{document}